\documentclass[letterpaper, 10 pt, conference]{ieeeconf}  

\IEEEoverridecommandlockouts                              
\usepackage{amsmath} 
\usepackage{amssymb}  
\usepackage{amsfonts}

\usepackage{graphicx}
\usepackage{epstopdf}
\usepackage{float} 
\usepackage[skip=0pt,font=footnotesize]{caption}

\usepackage[margin = 1in]{geometry}

\usepackage{subcaption}

\usepackage{etaremune}

\usepackage{tikz}
\usetikzlibrary{calc,patterns,decorations.pathmorphing,decorations.markings}
\usetikzlibrary{shapes,arrows}
\usepackage{verbatim}

\tikzstyle{block} = [draw, fill=blue!20, rectangle, 
    minimum height=3em, minimum width=6em]
\tikzstyle{sum} = [draw, fill=blue!20, circle, node distance=1cm]
\tikzstyle{input} = [coordinate]
\tikzstyle{output} = [coordinate]
\tikzstyle{pinstyle} = [pin edge={to-,thin,black}]

\usetikzlibrary{positioning}
\usetikzlibrary{math}

\usepackage{tikz}
\usetikzlibrary{shapes,arrows,calc,positioning}
\tikzstyle{bigblock} = [draw, fill=blue!20, rectangle, 
    minimum height=6em, minimum width=8em]
\tikzstyle{medblock} = [draw, fill=blue!20, rectangle, 
    minimum height=4em, minimum width=4em]    
\tikzstyle{mux} = [draw, fill=black!20, rectangle, 
    minimum height=5em, minimum width=0.1em]    
\tikzstyle{smallblock} = [draw, fill=blue!20, rectangle, 
    minimum height=3em, minimum width=4em]
\tikzstyle{sum} = [draw, fill=blue!20, circle, node distance=1cm]
\tikzstyle{signal} = [coordinate]
\tikzstyle{pinstyle} = [pin edge={to-,thin,black}]
\tikzstyle{block} = [draw, fill=blue!20, rectangle, 
    minimum height=3em, minimum width=6em]
\tikzstyle{blockS} = [draw, fill=blue!20, rectangle, 
    minimum height=3em, minimum width=4em]    
\tikzstyle{input} = [coordinate]
\tikzstyle{output} = [coordinate]
\usetikzlibrary{matrix}

\tikzset{add/.style n args={4}{
    minimum width=1mm,
    path picture={
        \draw[black, thick] 
            (path picture bounding box.south east) -- (path picture bounding box.north west)
            (path picture bounding box.south west) -- (path picture bounding box.north east);
        }
    }
}

\tikzset{Frame_into/.pic={
        code={\tikzset{scale=1}
        \tikzmath
            {
                \l  = 1;
                \Rc = 0.15;
            } 
        \draw [thick,->] (0, 0) -- +(\l, 0);
        \draw [thick,->] (0, 0) -- +(0, \l);
        \draw [thick, fill=white] 
			    (0,0) circle [radius=\Rc];
	    \draw [thick] ({\Rc*cos(45)}, {0\Rc*sin(45)}) 
	               -- ({\Rc*cos(225)}, {0\Rc*sin(225)});
        \draw [thick] ({\Rc*cos(135)}, {0\Rc*sin(135)}) 
	               -- ({\Rc*cos(315)}, {0\Rc*sin(315)});
  }}
}

\tikzset{Frame_outof/.pic={
        code={\tikzset{scale=1}
        \tikzmath
            {
                \l  = 1;
                \Rc = 0.15;
            } 
        \draw [thick,->] (0, 0) -- +(\l, 0);
        \draw [thick,->] (0, 0) -- +(0, \l);
        \draw [thick, fill=white] 
			    (0,0) circle [radius=\Rc];
	    \draw [thick, fill=black] 
			    (0,0) circle [radius={0.2*\Rc}];
  }}
}

\usepackage{hyperref}
\usepackage{xcolor}
\hypersetup{
    colorlinks,
    linkcolor={blue!100!black},
    citecolor={blue!50!black},
    urlcolor={blue!80!black}
}

\newcommand{\bc}{\begin{center}}
\newcommand{\ec}{\end{center}}
\newcommand{\benum}{\begin{enumerate}}
\newcommand{\eenum}{\end{enumerate}}
\newcommand{\nn}{\nonumber}
\newcommand{\matl}{\left[ \begin{array}}
\newcommand{\matr}{\end{array} \right]}
\newcommand{\matls}{\left[ \begin{smallmatrix}}
\newcommand{\matrs}{\end{smallmatrix} \right]}
\newcommand{\isdef}{\stackrel{\triangle}{=}}

\newcommand{\vect}[1]{\overset{\rightharpoonup}{#1}}

\newcommand{\rmA}{{\rm A}}
\newcommand{\rmB}{{\rm B}}

\newcommand{\rmT}{{\rm T}}

\newcommand{\rmc}{{\rm c}}
\newcommand{\rmd}{{\rm d}}

\newcommand{\rmf}{{\rm f}}

\newcommand{\BBR}{{\mathbb R}}

\newcommand{\SB}{{\mathcal B}}

\newcommand{\frameddot}[2]{\stackrel{{\rm #1}\bullet \bullet}{#2}}

\renewcommand{\matl}{\begin{bmatrix}}
\renewcommand{\matr}{\end{bmatrix} }

\newcommand{\tarrow}[1]{\overset{\rightarrow}{#1}}
\usepackage{graphics} 
\usepackage{amsmath} 
\usepackage{amssymb}  

\newcommand{\rotation}[2]{ {\substack{#1 \\ {\boldsymbol{\longrightarrow}} \\ #2}} }

\usepackage[style=ieee ,sorting=none,maxbibnames=99,giveninits, doi=false]{biblatex}
\title{\LARGE \bf
Longitudinal Flight Dynamics Control Based on Feedback Linearization and Normal Canonical Form
}

\title{\LARGE \bf
MIMO Input-Output Linearization with Applications for Longitudinal Flight Dynamics
}

\title{\LARGE \bf
Input-Output Linearization of MIMO Systems \\ with Application to Longitudinal Flight Dynamics
}

\title{Adaptive Finite-Time Control of Bicopter with Unknown Dynamics }

\title{Adaptive Nonlinear Control of a Bicopter with Unknown Dynamics }

\author{Jhon Manuel Portella Delgado and Ankit Goel
\thanks{Jhon Manuel Portella Delgado is a graduate student in the Department of Mechanical Engineering, University of Maryland, Baltimore County, 1000 Hilltop Circle, Baltimore, MD 21250. {\tt\small jportella@umbc.edu}}%
\thanks{Ankit Goel is an Assistant Professor in the Department of Mechanical Engineering, University of Maryland, Baltimore County,1000 Hilltop Circle, Baltimore, MD 21250. {\tt\small ankgoel@umbc.edu }}%
}
\usepackage{graphicx}      

\begin{document}

\maketitle

\begin{abstract}
%
This paper presents an adaptive, model-based, nonlinear controller for the bicopter trajectory-tracking problem. 
The nonlinear controller is constructed by dynamically extending the bicopter model, 
stabilizing the extended dynamics using input-output linearization,
augmenting the controller with a finite-time convergent parameter estimator, and designing a linear tracking controller. 
%
Unlike control systems based on the time separation principle to separate the translational and rotational dynamics, the proposed technique is applied to design a controller for the full nonlinear dynamics of the system to obtain the desired transient performance. 
The proposed controller is validated in simulation for a smooth and nonsmooth trajectory-tracking problem.

%
\end{abstract}
\textit{\bf keywords:} feedback linearization, dynamic extension,  finite-time estimation, multicopter.

\section{INTRODUCTION}




Multicopter UAVs have found great success as an inexpensive  tool in several engineering applications such as precision agriculture \cite{mukherjee2019}, environmental survey \cite{lucieer2014,klemas2015}, construction management \cite{li2019} and load transportation \cite{villa2020}. 
%
%
However, the low cost of building and operating multicopters fuels novel configurations and designs.
Furthermore, the intended operating envelope continually expands as the novel designs are used in novel applications.
Thus, due to nonlinear, time-varying, unmodeled dynamics, unknown and uncertain operating environments, and fast development cycles of novel configurations, multicopter control remains a challenging problem. 

Several control techniques have been applied to design control systems for multicopters \cite{nascimento2019,marshall2021,castillo2004}.
However, these techniques often require an accurate plant model and, thus, are susceptible to 
physical model parameter uncertainty \cite{emran2018,amin2016}.
Several adaptive control techniques have been applied to address the problem of unmodeled, unknown, and uncertain dynamics, such as model reference adaptive control \cite{whitehead2010,dydek2012}, L1 adaptive control \cite{zuo2014}, adaptive sliding mode control \cite{espinoza2021trajectory,wu20221, mofid2018}, retrospective cost adaptive control \cite{goel_adaptive_pid_2021,spencer2022}.

The dynamics of a multicopter consist of coupled translational and rotational dynamics, resulting in 12th-order nonlinear dynamics.  
The \textit{state-of-the-art} control architectures decompose the nonlinear dynamics into the linear translational and the nonlinear rotational dynamics \cite{px4_architecture}.
Cascaded stabilizing and regulating controllers are then designed for the translational and rotational dynamics separately in a multiloop architecture.  
%
%
%
%
Although controllers with theoretical performance guarantees can be designed and implemented for each loop in a cascaded control system, the performance of the closed-loop system can not be guaranteed. 
The cascaded control systems are based on the time separation principle to justify the cascaded loop architecture, which applies to the case where each successive loop is significantly faster than the previous loop. 
This crucial assumption allows the coupled dynamics to be decoupled and is used to design simpler controllers for each loop.  
However, as is well known, the entire control system fails if any loop fails. 
Since the controllers in each loop are often manually tuned, such cascaded control systems are difficult to tune and thus are highly susceptible to failure. 

In this paper, we thus consider the problem of designing an adaptive control system for the fully coupled nonlinear dynamics of a multicopter system. 
To simplify the presentation of the controller design technique, we consider a bicopter system, which retains the coupled nonlinear dynamics of a quadcopter system, but is modeled by a 6th-order nonlinear system instead of a 12th-order nonlinear system \cite{meinlschmidt2014cascaded,meinlschmidt2014cascaded2}.

The proposed controller is based on the input-output linearization (IOL) technique \cite{slotine1991applied, isidori1985nonlinear, khalil2015nonlinear}. 
To avoid the singularity of the resulting nonlinear input map by applying the IOL technique to the 6th-order nonlinear system, we first dynamically extend the system to design a linearizing controller \cite{triska2021dynamic}.  
Fortuitously, the dynamically extended system is fully linearizable and thus has no zero dynamics. 
Next, we design an adaptive parameter estimator with finite-time convergence property to rapidly estimate the uncertain parameters of the nonlinear dynamics.  
Finally, a static feedback controller is designed to obtain desired transient characteristics and track a desired trajectory. 

The contributions of this paper are thus 
1) the design of a linearizing controller for the fully nonlinear extended bicopter system without decoupling the nonlinear system into simpler subsystems, 
2) the adaptive extension of the proposed controller with a finite-time convergent parameter estimator and 
3) validation of the proposed controller in a smooth and nonsmooth trajectory-tracking problem.   

The paper is organized as follows. 
Section \ref{sec:prob_formulation} derives the equation of motion of the bicopter system.
Section \ref{sec:IOL} presents the adaptive input-output linearizing controller for the extended bicopter system.
Section \ref{sec:simulations} presents simulation results to validate the control system proposed in this work.
Finally, the paper concludes with a discussion of results and future research directions in section \ref{sec:conclusions}.

\section{Bicopter Dynamics}
\label{sec:prob_formulation}

\begin{figure}
    \centering
\tikzset{every picture/.style={line width=0.75pt}} 

\begin{tikzpicture}[x=0.75pt,y=0.75pt,yscale=-0.8,xscale=0.8]

\draw   (382.73,198.06) .. controls (384.44,196.15) and (387.38,195.98) .. (389.3,197.69) .. controls (391.22,199.41) and (391.38,202.35) .. (389.67,204.26) .. controls (387.96,206.18) and (385.02,206.35) .. (383.1,204.63) .. controls (381.18,202.92) and (381.02,199.98) .. (382.73,198.06) -- cycle ;
\draw    (385.37,137.84) -- (385.88,195.97) ;
\draw [shift={(385.35,135.84)}, rotate = 89.5] [fill={rgb, 255:red, 0; green, 0; blue, 0 }  ][line width=0.08]  [draw opacity=0] (12,-3) -- (0,0) -- (12,3) -- cycle    ;
\draw    (449.36,199.25) -- (391.23,200.02) ;
\draw [shift={(451.36,199.23)}, rotate = 179.24] [fill={rgb, 255:red, 0; green, 0; blue, 0 }  ][line width=0.08]  [draw opacity=0] (12,-3) -- (0,0) -- (12,3) -- cycle    ;
\draw  [fill={rgb, 255:red, 0; green, 0; blue, 0 }  ,fill opacity=1 ] (384.93,201.16) .. controls (384.93,200.46) and (385.5,199.9) .. (386.2,199.9) .. controls (386.9,199.9) and (387.47,200.46) .. (387.47,201.16) .. controls (387.47,201.86) and (386.9,202.43) .. (386.2,202.43) .. controls (385.5,202.43) and (384.93,201.86) .. (384.93,201.16) -- cycle ;
\draw    (322.35,83.75) -- (220.35,195.75) ;
\draw [shift={(271.35,139.75)}, rotate = 132.32] [color={rgb, 255:red, 0; green, 0; blue, 0 }  ][fill={rgb, 255:red, 0; green, 0; blue, 0 }  ][line width=0.75]      (0, 0) circle [x radius= 1.34, y radius= 1.34]   ;
\draw    (200.75,177.75) -- (220.35,195.75) ;
\draw   (183.73,197.05) .. controls (181.5,194.8) and (183.84,188.88) .. (188.95,183.81) .. controls (194.05,178.75) and (200,176.46) .. (202.22,178.71) .. controls (204.45,180.95) and (202.11,186.88) .. (197.01,191.94) .. controls (191.9,197.01) and (185.95,199.29) .. (183.73,197.05) -- cycle ;
\draw   (202.22,178.71) .. controls (200,176.46) and (202.33,170.54) .. (207.44,165.47) .. controls (212.55,160.41) and (218.49,158.12) .. (220.72,160.37) .. controls (222.94,162.61) and (220.61,168.54) .. (215.5,173.6) .. controls (210.39,178.66) and (204.45,180.95) .. (202.22,178.71) -- cycle ;
\draw    (302.75,65.75) -- (322.35,83.75) ;
\draw   (285.73,85.05) .. controls (283.5,82.8) and (285.84,76.88) .. (290.95,71.81) .. controls (296.05,66.75) and (302,64.46) .. (304.22,66.71) .. controls (306.45,68.95) and (304.11,74.88) .. (299.01,79.94) .. controls (293.9,85.01) and (287.95,87.29) .. (285.73,85.05) -- cycle ;
\draw   (304.22,66.71) .. controls (302,64.46) and (304.33,58.54) .. (309.44,53.47) .. controls (314.55,48.41) and (320.49,46.12) .. (322.72,48.37) .. controls (324.94,50.61) and (322.61,56.54) .. (317.5,61.6) .. controls (312.39,66.66) and (306.45,68.95) .. (304.22,66.71) -- cycle ;
\draw  [dash pattern={on 4.5pt off 4.5pt}]  (331.49,138.96) -- (271.35,139.75) ;
\draw  [fill={rgb, 255:red, 0; green, 0; blue, 0 }  ,fill opacity=1 ] (326.07,201.55) .. controls (326.07,200.98) and (326.53,200.51) .. (327.1,200.51) .. controls (327.68,200.51) and (328.14,200.98) .. (328.14,201.55) .. controls (328.14,202.12) and (327.68,202.59) .. (327.1,202.59) .. controls (326.53,202.59) and (326.07,202.12) .. (326.07,201.55) -- cycle ;
\draw [line width=1.5]    (410.42,70.59) -- (410.12,109.92) ;
\draw [shift={(410.09,113.92)}, rotate = 270.44] [fill={rgb, 255:red, 0; green, 0; blue, 0 }  ][line width=0.08]  [draw opacity=0] (8.75,-4.2) -- (0,0) -- (8.75,4.2) -- (5.81,0) -- cycle    ;
\draw [line width=1.5]    (196.22,172.04) -- (172.07,149.28) ;
\draw [shift={(169.16,146.54)}, rotate = 43.29] [fill={rgb, 255:red, 0; green, 0; blue, 0 }  ][line width=0.08]  [draw opacity=0] (8.75,-4.2) -- (0,0) -- (8.75,4.2) -- (5.81,0) -- cycle    ;
\draw [line width=1.5]    (298.42,60.42) -- (274.27,37.66) ;
\draw [shift={(271.35,34.92)}, rotate = 43.29] [fill={rgb, 255:red, 0; green, 0; blue, 0 }  ][line width=0.08]  [draw opacity=0] (8.75,-4.2) -- (0,0) -- (8.75,4.2) -- (5.81,0) -- cycle    ;
\draw   (192.75,92.55) .. controls (192.54,89.99) and (194.46,87.75) .. (197.02,87.55) .. controls (199.58,87.35) and (201.82,89.26) .. (202.02,91.83) .. controls (202.22,94.39) and (200.31,96.63) .. (197.75,96.83) .. controls (195.18,97.03) and (192.95,95.11) .. (192.75,92.55) -- cycle ;
\draw    (151.07,49) -- (193.41,88.83) ;
\draw [shift={(149.61,47.63)}, rotate = 43.25] [fill={rgb, 255:red, 0; green, 0; blue, 0 }  ][line width=0.08]  [draw opacity=0] (12,-3) -- (0,0) -- (12,3) -- cycle    ;
\draw    (239.68,45.24) -- (200.04,87.77) ;
\draw [shift={(241.05,43.78)}, rotate = 132.99] [fill={rgb, 255:red, 0; green, 0; blue, 0 }  ][line width=0.08]  [draw opacity=0] (12,-3) -- (0,0) -- (12,3) -- cycle    ;
\draw  [fill={rgb, 255:red, 0; green, 0; blue, 0 }  ,fill opacity=1 ] (196.51,93.1) .. controls (196,92.62) and (195.99,91.82) .. (196.47,91.31) .. controls (196.95,90.81) and (197.75,90.79) .. (198.26,91.27) .. controls (198.77,91.76) and (198.78,92.56) .. (198.3,93.07) .. controls (197.82,93.57) and (197.01,93.59) .. (196.51,93.1) -- cycle ;
\draw  [draw opacity=0] (295.51,113.84) .. controls (295.56,113.84) and (295.62,113.84) .. (295.67,113.84) .. controls (301.93,113.99) and (306.85,120.97) .. (306.66,129.45) .. controls (306.57,133.33) and (305.43,136.85) .. (303.61,139.52) -- (295.32,129.19) -- cycle ; \draw    (295.67,113.84) .. controls (301.93,113.99) and (306.85,120.97) .. (306.66,129.45) .. controls (306.59,132.56) and (305.85,135.43) .. (304.62,137.82) ; \draw [shift={(303.61,139.52)}, rotate = 291.91] [fill={rgb, 255:red, 0; green, 0; blue, 0 }  ][line width=0.08]  [draw opacity=0] (7.2,-1.8) -- (0,0) -- (7.2,1.8) -- cycle    ; \draw [shift={(295.51,113.84)}, rotate = 24.68] [fill={rgb, 255:red, 0; green, 0; blue, 0 }  ][line width=0.08]  [draw opacity=0] (7.2,-1.8) -- (0,0) -- (7.2,1.8) -- cycle    ;

\draw (413.67,79.33) node [anchor=north west][inner sep=0.75pt]   [align=left] {{\small $\vect g$}};
\draw (453,190) node [anchor=north west][inner sep=0.75pt]   [align=left] {{\small $\hat{i}_{\rm{A}}$}};
\draw (375,120) node [anchor=north west][inner sep=0.75pt]   [align=left] {{\small $\hat{j}_{\rm{A}}$}};
\draw (234,31) node [anchor=north west][inner sep=0.75pt]  [rotate=-313.75] [align=left] {{\small $\hat{i}_{\rm B}$}};
\draw (129.9,40) node [anchor=north west][inner sep=0.75pt]  [rotate=-313.75] [align=left] {{\small $\hat{j}_{\rm B}$}};
\draw (287,20.33) node [anchor=north west][inner sep=0.75pt]   [align=left] {{\small ${\vect f}_2$}};
\draw (184,135.33) node [anchor=north west][inner sep=0.75pt]   [align=left] {{\small ${\vect f}_1$}};
\draw (260.67,125) node [anchor=north west][inner sep=0.75pt]   [align=left] {{\small c}};
\draw (327.33,200) node [anchor=north west][inner sep=0.75pt]   [align=left] {{\small \textit{w}}};
\draw (308.33,115) node [anchor=north west][inner sep=0.75pt]   [align=left] {{\small $\theta$}};
\end{tikzpicture}
\caption{Bicopter configuration considered in this paper. The bicopter is constrained to the $\hat \imath _\rmA-\hat \jmath _\rmA$ plane and rotates about the $\hat k_\rmA$ axis of the inertial frame $\rm F_A.$  }
    \label{fig:Bicopter}
    \vspace{-2em}
\end{figure}
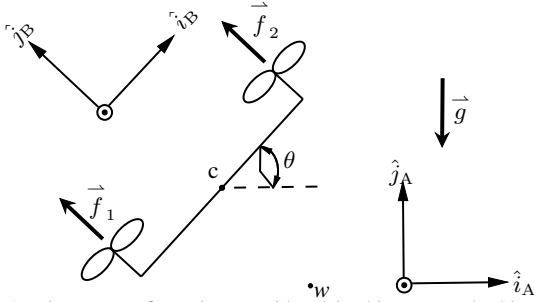

%
This section derives the equation of motion of a bicopter using Newton-Euler dynamics.
Let $\rm F_A$ be an inertial frame and let $\rm F_B$ be a frame fixed to the bicopter $\SB$ as shown in Figure \ref{fig:Bicopter}. 
Note that $\rm F_B$ is obtained by rotating it about the $\hat k_\rmA$ axis of $\rm F_A$ by $\theta,$ and thus 
\begin{align}
    {\rm F_A}
        \rotation{\theta}{3}
    {\rm F_B}.
\end{align}

Letting $\rmc$ denote the center of mass of the bicopter and $w$ denote a fixed point on Earth, it follows from Newton's second law that 
\begin{align}
    m
    \frameddot{A}{\vect r}_{\rmc/w} 
        &=
            m \vect g + \vect f,
    \label{eq:N2L_pos}
\end{align}
where 
$m$ is the mass of the bicopter,
$\vect g$ is the acceleration due to gravity, 
and $\vect f$ is the total force applied by the propellers to the bicopter. 
Letting $\vect f_1 = f_1 \hat \jmath_\rmB$ and $\vect f_2 = f_2 \hat \jmath_\rmB$ denote the forces applied by the two propellers, it follows that $\vect f =  f_1 \hat j_{\rm B} + f_2 \hat j_\rmB $.
Writing
$\vect r_{\rmc/w} 
        =
            r_1 \hat \imath_\rmA +
            r_2 \hat \jmath_\rmA$ yields
\begin{align}
    m \ddot r_1 &= -(f_1 + f_2) \sin \theta, 
    \label{eq:eom_r1}
    \\
    m \ddot r_2 &= (f_1 + f_2) \cos \theta - m g.
    \label{eq:eom_r2}
\end{align}
Next, it follows from Euler's equation that 
\begin{align}
    \tarrow J_{\SB/\rmc} \frameddot{A}{\vect \omega}_{\rm B/A} 
        =
            \vect M_{\SB/\rmc}.
    \label{eq:EulersEqn}
\end{align}
Note that  $\tarrow J_{\SB/\rmc} \frameddot{A}{\vect \omega}_{\rm B/A}  = J \ddot \theta \hat k_\rmB $
and $\vect M_{\SB/\rmc} = \ell (f_2 - f_1) \hat k_\rmB,$ where $\ell$ is the length of the bicopter arm, 
and thus it follows from \eqref{eq:EulersEqn} that
\begin{align}
    J \ddot \theta = \ell (f_2 - f_1). 
    \label{eq:eom_theta}
\end{align}

The equations of motion of the bicopter, given by \eqref{eq:eom_r1}, \eqref{eq:eom_r2}, and \eqref{eq:eom_theta}, can be written in state-space form as 
\begin{align}
    \dot x = f(x) + g(x)u,
    \label{eq:former_ss}
\end{align}
where 
\begin{align}
    x 
    &\isdef 
        \matl 
            r_1 &
            r_2 &
            \theta &
            \dot r_1 &
            \dot r_2 &
            \dot \theta
        \matr^\rmT,
    \label{initial_states}\\
    u
    &\isdef 
        \matl 
            f_1 + f_2 &
            \ell (f_2 - f_1)
        \matr^\rmT 
        \label{control_signals}
        ,
\end{align}
and 
\begin{align}
    f(x)
        \isdef
            \matl
                x_4\\
                x_5\\
                x_6 \\
                0\\
                -g\\
                0
            \matr, \quad 
    g(x)
        \isdef
            \matl
                0 & 0\\
                0 & 0\\
                0 & 0\\
                \dfrac{-\sin{x_3}}{m} & 0 \\[0.2cm]
                \dfrac{\cos{x_3}}{m} & 0\\
                0 & \dfrac{1}{J}
            \matr.
\end{align}

\section{Adaptive Dynamic Input-Output Linearizing Control}
\label{sec:IOL}

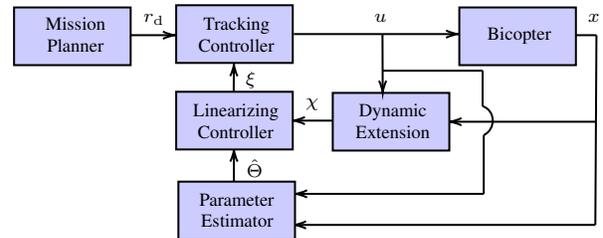
\begin{figure}[!ht] 

\tikzset{every picture/.style={line width=0.75pt}} 

\begin{tikzpicture}[x=0.75pt,y=0.75pt,yscale=-0.56,xscale=0.56]
\tikzstyle{every node}=[font=\scriptsize]

\draw    (118.2,53.6) -- (156.2,53.6) ;
\draw [shift={(158.2,53.6)}, rotate = 180] [color={rgb, 255:red, 0; green, 0; blue, 0 }  ][line width=0.75]    (10.93,-3.29) .. controls (6.95,-1.4) and (3.31,-0.3) .. (0,0) .. controls (3.31,0.3) and (6.95,1.4) .. (10.93,3.29)   ;
\draw    (265.2,52.6) -- (411.2,52.6) ;
\draw [shift={(413.2,52.6)}, rotate = 180] [color={rgb, 255:red, 0; green, 0; blue, 0 }  ][line width=0.75]    (10.93,-3.29) .. controls (6.95,-1.4) and (3.31,-0.3) .. (0,0) .. controls (3.31,0.3) and (6.95,1.4) .. (10.93,3.29)   ;
\draw    (537.7,53.6) -- (536.7,224.6) ;
\draw    (537.2,131.6) -- (407.2,131.6) ;
\draw [shift={(405.2,131.6)}, rotate = 360] [color={rgb, 255:red, 0; green, 0; blue, 0 }  ][line width=0.75]    (10.93,-3.29) .. controls (6.95,-1.4) and (3.31,-0.3) .. (0,0) .. controls (3.31,0.3) and (6.95,1.4) .. (10.93,3.29)   ;
\draw    (435.2,117.6) .. controls (447.2,121.6) and (447.2,138.6) .. (436.2,143.6) ;
\draw    (436.2,143.6) -- (436.2,196.6) ;
\draw    (436.2,196.6) -- (268.2,196.6) ;
\draw [shift={(266.2,196.6)}, rotate = 360] [color={rgb, 255:red, 0; green, 0; blue, 0 }  ][line width=0.75]    (10.93,-3.29) .. controls (6.95,-1.4) and (3.31,-0.3) .. (0,0) .. controls (3.31,0.3) and (6.95,1.4) .. (10.93,3.29)   ;
\draw    (435.2,85.6) -- (435.2,117.6) ;
\draw    (536.7,224.6) -- (267.2,224.6) ;
\draw [shift={(265.2,224.6)}, rotate = 360] [color={rgb, 255:red, 0; green, 0; blue, 0 }  ][line width=0.75]    (10.93,-3.29) .. controls (6.95,-1.4) and (3.31,-0.3) .. (0,0) .. controls (3.31,0.3) and (6.95,1.4) .. (10.93,3.29)   ;
\draw    (211.2,185.6) -- (211.2,158.6) ;
\draw [shift={(211.2,156.6)}, rotate = 90] [color={rgb, 255:red, 0; green, 0; blue, 0 }  ][line width=0.75]    (10.93,-3.29) .. controls (6.95,-1.4) and (3.31,-0.3) .. (0,0) .. controls (3.31,0.3) and (6.95,1.4) .. (10.93,3.29)   ;
\draw    (211,105) -- (211.18,81.6) ;
\draw [shift={(211.2,79.6)}, rotate = 90.45] [color={rgb, 255:red, 0; green, 0; blue, 0 }  ][line width=0.75]    (10.93,-3.29) .. controls (6.95,-1.4) and (3.31,-0.3) .. (0,0) .. controls (3.31,0.3) and (6.95,1.4) .. (10.93,3.29)   ;
\draw    (300.2,130.6) -- (266.2,130.6) ;
\draw [shift={(264.2,130.6)}, rotate = 360] [color={rgb, 255:red, 0; green, 0; blue, 0 }  ][line width=0.75]    (10.93,-3.29) .. controls (6.95,-1.4) and (3.31,-0.3) .. (0,0) .. controls (3.31,0.3) and (6.95,1.4) .. (10.93,3.29)   ;
\draw  [fill=blue!20  ,fill opacity=1 ] (300.2,105.6) -- (405.2,105.6) -- (405.2,157.6) -- (300.2,157.6) -- cycle ;
\draw  [fill=blue!20  ,fill opacity=1 ] (159.2,26.6) -- (264.2,26.6) -- (264.2,78.6) -- (159.2,78.6) -- cycle ;
\draw  [fill=blue!20  ,fill opacity=1 ] (159.2,104.6) -- (264.2,104.6) -- (264.2,156.6) -- (159.2,156.6) -- cycle ;
\draw  [fill=blue!20  ,fill opacity=1 ] (161.2,185.6) -- (266.2,185.6) -- (266.2,237.6) -- (161.2,237.6) -- cycle ;
\draw  [fill=blue!20  ,fill opacity=1 ] (13.2,28.6) -- (118.2,28.6) -- (118.2,80.6) -- (13.2,80.6) -- cycle ;
\draw  [fill=blue!20  ,fill opacity=1 ] (414.2,28.6) -- (519.2,28.6) -- (519.2,80.6) -- (414.2,80.6) -- cycle ;
\draw    (519.2,53.6) -- (537.7,53.6) ;
\draw    (345,53) -- (345.19,103.6) ;
\draw [shift={(345.2,105.6)}, rotate = 269.78] [color={rgb, 255:red, 0; green, 0; blue, 0 }  ][line width=0.75]    (10.93,-3.29) .. controls (6.95,-1.4) and (3.31,-0.3) .. (0,0) .. controls (3.31,0.3) and (6.95,1.4) .. (10.93,3.29)   ;
\draw    (345.2,85.6) -- (435.2,85.6) ;

\draw (18,35) node [anchor=north west][inner sep=0.75pt]   [align=left] {\begin{minipage}[lt]{40.13pt}\setlength\topsep{0pt}
\begin{center}
Mission \\Planner
\end{center}

\end{minipage}};
\draw (145,113) node [anchor=north west][inner sep=0.75pt]   [align=left] {\begin{minipage}[lt]{54.89pt}\setlength\topsep{0pt}
\begin{center}
Linearizing \\Controller
\end{center}

\end{minipage}};
\draw (155,35) node [anchor=north west][inner sep=0.75pt]   [align=left] {\begin{minipage}[lt]{46.94pt}\setlength\topsep{0pt}
\begin{center}
Tracking\\Controller
\end{center}

\end{minipage}};
\draw (295,113) node [anchor=north west][inner sep=0.75pt]   [align=left] {\begin{minipage}[lt]{47.52pt}\setlength\topsep{0pt}
\begin{center}
Dynamic \\Extension
\end{center}

\end{minipage}};
\draw (148,193) node [anchor=north west][inner sep=0.75pt]   [align=left] {\begin{minipage}[lt]{53.18pt}\setlength\topsep{0pt}
\begin{center}
 Parameter\\Estimator
\end{center}

\end{minipage}};
\draw (417.6,43.6) node [anchor=north west][inner sep=0.75pt]   [align=left] {\begin{minipage}[lt]{40.15pt}\setlength\topsep{0pt}
\begin{center}
Bicopter
\end{center}

\end{minipage}};
\draw (335,32) node [anchor=north west][inner sep=0.75pt]   [align=left] {$u$};
\draw (220,163) node [anchor=north west][inner sep=0.75pt]   [align=left] {$\hat{\Theta}$};
\draw (273,108) node [anchor=north west][inner sep=0.75pt]   [align=left] {$\chi$};
\draw (219,84) node [anchor=north west][inner sep=0.75pt]   [align=left] {$\xi$};
\draw (127,32) node [anchor=north west][inner sep=0.75pt]   [align=left] {$r_{\rm d}$};
\draw (527.2,32) node [anchor=north west][inner sep=0.75pt]   [align=left] {$x$};

\end{tikzpicture}

\vspace{1pt}
\caption{Adaptive, Dynamic, Input-Output Linearizing Control Architecture.}
    \label{fig:ADIOL_architecture}
\end{figure}

This section develops the adaptive, dynamic, input-output linearizing (ADIOL) controller for a bicopter with unknown dynamics to track a desired trajectory. 
The ADIOL controller is designed by dynamically extending the bicopter system, designing an input-output linearizing controller for the extended system, and then augmenting the controller with a finite-time convergent parameter estimator.
In the bicopter system, dynamic extension is necessary to avoid the singularity of the input map that needs inversion to design the linearizing controller. 
Finally, a linear controller is designed for the linearized system to track a desired trajectory. 
Figure \ref{fig:ADIOL_architecture} shows the architecture of the ADIOL control system.

%
The input-output linearization technique is based on the system's normal form \cite[p.~517]{Khalil:1173048}, \cite{brunovsky1970classification,fliess1993some, vcelikovsky1995global}.  
In the normal form, the system's dynamics are decomposed into coupled zero dynamics and linearizable dynamics. 
%
%
However, successful linearization of the linearizable dynamics requires the inversion of a nonlinear map, which may or may not be possible. 
Moreover, the zero dynamics must be asymptotically stable for the input-output linearizing controller to be useful. 
In the following, we use the notation presented in Section III of \cite{portella2024circumventing} to design an input-output linearizing controller.

%


%

%

In the bicopter dynamics, unfortunately, the input-output linearizing controller can not be constructed due to the non-invertibility of the nonlinear map $\beta(x),$ as shown below. 
%
%
%
Let the position of the bicopter be the output of the system, that is, $y=h(x) = \matl r_1 & r_2\matr^{\rm T}.$
Note that the relative degrees $\rho_1$ and $\rho_2$ of $y_1$ and $y_2,$ respectively, are $2$, and thus
\begin{align}
    \beta(x)
        =
    \matl
        L_g L_f h_1^{\rho_1-1}(x)\\
        L_g L_f h_2^{\rho_2-1}(x)
    \matr 
     =
     \matl
        -\dfrac{\sin(x_3)}{m} & 0\\[2mm]
        \dfrac{\cos(x_3)}{m} & 0
     \matr.
\end{align}
Since $\beta(x)$ is non-invertible, the bicopter system can not be input-output linearized with the current choice of inputs and outputs.

\subsection{Dynamic Extension}
\label{sec:DIOL}

To overcome this restriction caused by the singularity of $\beta$,
%
we dynamically extend the bicopter system as shown below. 
By defining 
%
\begin{align}
    w 
        \isdef 
            \matl 
                \ddot u_1 \\
                u_2
            \matr,
    \quad 
    \chi 
        \isdef  
            \matl
                x \\
                u_1 \\
                \dot u_1 
            \matr,
    \label{dynamic_extension_definition}
\end{align}
the dynamically extended system is
\begin{align}
    \dot \chi = F(\chi) + G(\chi) w,
    \label{eq:extended_system}
\end{align}
where 
\begin{align}
    F(\chi) 
        \isdef  
            \matl 
                \chi_4\\
                \chi_5\\
                \chi_6\\
                -\dfrac{\sin{\chi_3}}{m}\chi_7\\[0.2cm]
                -g+\dfrac{\cos{\chi_3}}{m}\chi_7\\
                0\\
                \chi_8\\
                0
            \matr
    ,
    G(\chi) 
        \matl 
            0 & 0\\
            0 & 0\\
            0 & 0\\
            0 & 0\\
            0 & 0\\
            0 & \dfrac{1}{J}\\
            0 & 0\\
            1 & 0
        \matr.    
\end{align}

Let the position of the bicopter be the output of the dynamically extended system \eqref{eq:extended_system}, that is, $y=H(\chi) = \matl r_1 & r_2\matr^{\rm T}.$
Note that the relative degrees $\rho_1$ and $\rho_2$ of $y_1$ and $y_2$ are $4$, respectively, and thus 
\begin{align}
    \beta(\chi)
        &\isdef 
            \matl
                L_{G}L_{F}^{3}H_1(\chi)\\
                L_{G}L_{F}^{3}H_2(\chi)
            \matr
        =
            \matl
                -\dfrac{\sin{\chi_3}}{m} &  -\dfrac{\cos{\chi_3}}{m}\dfrac{\chi_7}{J}\\
                \dfrac{\cos{\chi_3}}{m} & -\dfrac{\sin{\chi_3}}{m}\dfrac{\chi_7}{J}
            \matr,
\end{align}
Assuming $\chi_7 = u_1 \neq 0,$ which is a reasonable assumption in multicopter control, 
\begin{align}
    \beta^{-1}(\chi)
        =
            \matl
                -m\,\sin \left(\chi_3 \right) & m\,\cos \left(\chi_3 \right)\\
                -\dfrac{J\,m\,\cos \left(\chi_3 \right)}{\chi_7 } & -\dfrac{J\,m\,\sin \left(\chi_3 \right)}{\chi_7 }
            \matr.
\end{align}
Furthermore, the relative degree $\rho$ of the output in the extended system is $ 8,$ which  is equal to the dimension of the extended system state $\chi$. 
Therefore, there are no zero dynamics in the extended bicopter dynamics.

Finally, the normal form of \eqref{eq:extended_system} is 
\begin{align}
    \dot \xi
        = 
            A_\rmc \xi + B_\rmc (\alpha(\chi) + \beta(\chi) u ),\label{dot_xi_formulation}
\end{align}
where 
\begin{align}
    \xi 
        \isdef 
            \matl
                H_1(\chi)\\
                L_F H_1(\chi) \\
                \vdots\\
                L_F H_1^{\rho_1-1}(\chi)
                \\
                H_2(\chi)\\
                L_F H_2(\chi) \\
                \vdots\\
                L_F H_2^{\rho_2-1}(\chi)
            \matr
        =
            \matl
                \chi_1\\
                \chi_4\\
                -\dfrac{\sin{(\chi_3)}}{m}\chi_7\\
                -\dfrac{\cos{(\chi_3)}}{m}\chi_6\chi_7 -\dfrac{\sin{(\chi_3)}}{m}\chi_8
                \\
                \chi_2\\
                \chi_5\\
                -g+\dfrac{\cos{(\chi_3)}}{m}\chi_7\\
                -\dfrac{\sin{(\chi_3)}}{m}\chi_6\chi_7 + \dfrac{\cos{\chi_3}}{m}\chi_8
            \matr,
    \label{xi_definition}
\end{align}
\begin{align}
    A_\rmc  \isdef 
        \matl
            0 & 1 & 0 & 0 & 0 & 0 & 0 & 0\\
            0 & 0 & 1 & 0 & 0 & 0 & 0 & 0\\
            0 & 0 & 0 & 1 & 0 & 0 & 0 & 0\\
            0 & 0 & 0 & 0 & 0 & 0 & 0 & 0\\
            0 & 0 & 0 & 0 & 0 & 1 & 0 & 0\\
            0 & 0 & 0 & 0 & 0 & 0 & 1 & 0\\
            0 & 0 & 0 & 0 & 0 & 0 & 0 & 1\\
            0 & 0 & 0 & 0 & 0 & 0 & 0 & 0
        \matr,
    B_\rmc  \isdef 
        \matl
         0 & 0\\
         0 & 0\\
         0 & 0\\
         1 & 0\\
         0 & 0\\
         0 & 0\\
         0 & 0\\
         0 & 1\\
        \matr,
        \label{A_and_B_normal_form}
\end{align}
and 
\begin{align}
    \alpha(\chi)
        &\isdef
            \matl
                L_{F}^4H_1(\chi)\\
                L_{F}^4H_2(\chi)
            \matr\nonumber\\
    &=
        \matl       
            -\dfrac{\chi_6 \,{\left(2\,\chi_8 \,\cos \left(\chi_3 \right)-\chi_6 \,\chi_7 \,\sin \left(x_3 \right)\right)}}{m}\\
            -\dfrac{\chi_6 \,{\left(2\,\chi_8 \,\sin \left(\chi_3 \right)+\chi_6 \,\chi_7 \,\cos \left(\chi_3 \right)\right)}}{m}
         \matr.
\end{align}

Since the map $\beta$ is invertible, the control law 
\begin{align}
    w = -\beta(\chi)^{-1}(\alpha(\chi) - v),
    \label{eq:IOL_control_law}
\end{align}
yields the linearized dynamics
\begin{align}
    \dot \xi = A_\rmc  \xi + B_\rmc  v.
    \label{eq:linearized_dynamics}
\end{align}

\subsection{Adaptive Augmentation with Finite-Time Convergence}

The linearizing control law \eqref{eq:IOL_control_law} requires the precise knowledge of the mass $m$ and the inertia $J$ of the bicopter.
However, since these physical parameters are not typically precisely known or are often time-varying, we use an online parameter estimator with finite-time convergence property to estimate the unknown parameters.
The finite-time estimation of the parameters thus allows the exact linearization of the extended bicopter dynamics. 

The finite-time convergent parameter estimator \eqref{eq:estimator} is constructed by modifying the gradient-based exponentially stable parameter estimator using the finite-time optimization theory presented in \cite{garg2020fixed,garg2020cappa}. 

To estimate the unknown parameters, we use the bicopter system \eqref{eq:former_ss} to formulate a linear regressor equation, as shown below. 
This approach is motivated by the DREM parameter estimator described in \cite{aranovskiy2016parameters}.
In particular, we write \eqref{eq:former_ss} as
\begin{align}
    \dot x - \Psi(x) = \Phi(x,u) \Theta, 
    \label{eq:extended_system_linear_para}
\end{align}
where 
\begin{align}
    \Psi(x)
        \isdef 
            \matl 
                x_4\\
                x_5\\
                x_6\\
                0\\
                -g\\
                0
            \matr, 
    \Phi(x,u)
        \isdef 
            \matl 
                0 & 0\\
                0 & 0\\
                0 & 0\\
                -\sin{(x_3)}u_1 & 0\\
                \cos{(x_3)}u_1 & 0\\
                0 & u_2
            \matr,             
\end{align}
and 
\begin{align}
    \Theta 
        \isdef  
            \matl 
                m^{-1} \\
                J^{-1}
            \matr.
\end{align}
Since $x$ and $u$ are assumed to be known, the signals $\Psi$ and $\Phi$ can be computed online.
However, $\dot x$ is unknown. 
To avoid the requirement of $\dot x$, we filter \eqref{eq:extended_system_linear_para} with a strictly proper filter $R(s)$ to obtain the linear regressor
\begin{align}
    x_\rmf = \Phi_\rmf \Theta,
    \label{eq:LRE}
\end{align}
where 
\begin{align}
    x_\rmf
        &\isdef
        R(s) [\dot x - \Psi(x)]
    , \quad
    \Phi_\rmf 
        \isdef R(s) [\Phi(x,u)].
\end{align}
For example, letting $R(s) = \dfrac{1}{s+\gamma},$ where $\gamma>0,$ yields
\begin{align}
    x_\rmf &= \dfrac{sx}{s+\gamma} - \dfrac{\Psi(x)}{s+\gamma} , \quad 
    \Phi_\rmf = \dfrac{\Phi(x,u)}{s+\gamma},
\end{align}
where $x_\rmf$ and $\Phi_\rmf$ can now be computed online using only the measurements of the state $x$ and the input $u.$

Finally, we use the estimator
\begin{align}
    \dot{\hat{\Theta}}
        =
            -c_1 \dfrac{\Xi}{||\Xi||_2^{1-\alpha_1}}
            -c_2 \dfrac{\Xi}{||\Xi||_2^{1-\alpha_2}},
    \label{eq:estimator}
\end{align}
where 
\begin{align}
    \Xi
        \isdef 
            \overline{\Phi} \hat{\Theta} -\overline{X},
\end{align}
$c_1,c_2>0$, $0<\alpha_1<1$, and $\alpha_2>1,$ and 
\begin{align}
    \dot{\overline{X}} 
        &=
            -\lambda \overline{X} + 
            \Phi_\rmf^\rmT x_\rmf
        ,\\
    \dot{\overline{\Phi}} 
        &=
            -\lambda \overline{\Phi} + 
            \Phi_\rmf^{\rm T}\Phi_\rmf.
\end{align}
The positive scalar $\lambda>0$ is the exponential forgetting factor.
%
%
The data matrices $\overline{X} \in \BBR^{2} $ and $\overline{\Phi} \in \BBR^{2 \times 2}$ are initialized at zero.

The adaptive control law is thus \eqref{eq:IOL_control_law}, where $\alpha$ and $\beta$ are computed with the estimates of $m$ and $J$ given by \eqref{eq:estimator}, and $v$ is given by a linear tracking controller as described in the next section.

\subsection{Trajectory Tracking Control}



To track a desired trajectory, we use the following control laws. 

\subsubsection{Full-state feedback with Integral Action}
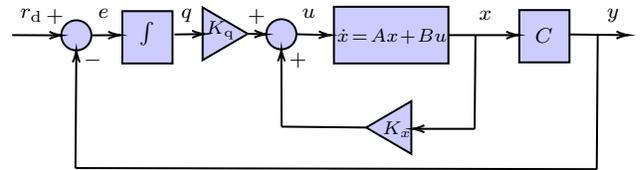
\begin{figure}[H]
    \centering
\tikzset{every picture/.style={line width=0.75pt}} 

\begin{tikzpicture}[x=0.75pt,y=0.75pt,yscale=-0.55,xscale=0.55]
\tikzstyle{every node}=[font=\footnotesize]

\draw  [fill=blue!20  ,fill opacity=1 ] (344.2,35.6) -- (449.2,35.6) -- (449.2,87.6) -- (344.2,87.6) -- cycle ;
\draw  [fill=blue!20  ,fill opacity=1 ] (514.2,38.6) -- (560.2,38.6) -- (560.2,86.6) -- (514.2,86.6) -- cycle ;
\draw    (449.2,61.6) -- (511.2,61.6) ;
\draw [shift={(513.2,61.6)}, rotate = 180] [color={rgb, 255:red, 0; green, 0; blue, 0 }  ][line width=0.75]    (10.93,-3.29) .. controls (6.95,-1.4) and (3.31,-0.3) .. (0,0) .. controls (3.31,0.3) and (6.95,1.4) .. (10.93,3.29)   ;
\draw  [fill=blue!20  ,fill opacity=1 ] (373.8,146.4) -- (415.2,122.2) -- (415.2,170.6) -- cycle ;
\draw    (474,62) -- (473.2,147.6) ;
\draw    (473.2,147.6) -- (418.2,147.6) ;
\draw [shift={(416.2,147.6)}, rotate = 360] [color={rgb, 255:red, 0; green, 0; blue, 0 }  ][line width=0.75]    (10.93,-3.29) .. controls (6.95,-1.4) and (3.31,-0.3) .. (0,0) .. controls (3.31,0.3) and (6.95,1.4) .. (10.93,3.29)   ;
\draw  [fill=blue!20  ,fill opacity=1 ] (282.2,60.8) .. controls (282.2,53.18) and (288.38,47) .. (296,47) .. controls (303.62,47) and (309.8,53.18) .. (309.8,60.8) .. controls (309.8,68.42) and (303.62,74.6) .. (296,74.6) .. controls (288.38,74.6) and (282.2,68.42) .. (282.2,60.8) -- cycle ;
\draw    (309.8,60.8) -- (339.2,61.55) ;
\draw [shift={(341.2,61.6)}, rotate = 181.46] [color={rgb, 255:red, 0; green, 0; blue, 0 }  ][line width=0.75]    (10.93,-3.29) .. controls (6.95,-1.4) and (3.31,-0.3) .. (0,0) .. controls (3.31,0.3) and (6.95,1.4) .. (10.93,3.29)   ;
\draw    (296.2,145.6) -- (373.8,146.4) ;
\draw    (296.2,145.6) -- (296.01,76.6) ;
\draw [shift={(296,74.6)}, rotate = 89.84] [color={rgb, 255:red, 0; green, 0; blue, 0 }  ][line width=0.75]    (10.93,-3.29) .. controls (6.95,-1.4) and (3.31,-0.3) .. (0,0) .. controls (3.31,0.3) and (6.95,1.4) .. (10.93,3.29)   ;
\draw  [fill=blue!20  ,fill opacity=1 ] (94.2,61.8) .. controls (94.2,54.18) and (100.38,48) .. (108,48) .. controls (115.62,48) and (121.8,54.18) .. (121.8,61.8) .. controls (121.8,69.42) and (115.62,75.6) .. (108,75.6) .. controls (100.38,75.6) and (94.2,69.42) .. (94.2,61.8) -- cycle ;
\draw  [fill=blue!20  ,fill opacity=1 ] (150.2,39.6) -- (196.2,39.6) -- (196.2,87.6) -- (150.2,87.6) -- cycle ;
\draw  [fill=blue!20  ,fill opacity=1 ] (265.2,60.4) -- (223.8,84.6) -- (223.8,36.2) -- cycle ;
\draw    (560,61) -- (612.2,61.58) ;
\draw [shift={(614.2,61.6)}, rotate = 180.61] [color={rgb, 255:red, 0; green, 0; blue, 0 }  ][line width=0.75]    (10.93,-3.29) .. controls (6.95,-1.4) and (3.31,-0.3) .. (0,0) .. controls (3.31,0.3) and (6.95,1.4) .. (10.93,3.29)   ;
\draw    (586.1,61.3) -- (585.2,184.6) ;
\draw    (107.2,183.6) -- (585.2,184.6) ;
\draw    (107.2,183.6) -- (107.99,77.6) ;
\draw [shift={(108,75.6)}, rotate = 90.42] [color={rgb, 255:red, 0; green, 0; blue, 0 }  ][line width=0.75]    (10.93,-3.29) .. controls (6.95,-1.4) and (3.31,-0.3) .. (0,0) .. controls (3.31,0.3) and (6.95,1.4) .. (10.93,3.29)   ;
\draw    (121.8,61.8) -- (146.2,61.62) ;
\draw [shift={(148.2,61.6)}, rotate = 179.57] [color={rgb, 255:red, 0; green, 0; blue, 0 }  ][line width=0.75]    (10.93,-3.29) .. controls (6.95,-1.4) and (3.31,-0.3) .. (0,0) .. controls (3.31,0.3) and (6.95,1.4) .. (10.93,3.29)   ;
\draw    (198,61) -- (223.2,60.63) ;
\draw [shift={(225.2,60.6)}, rotate = 179.16] [color={rgb, 255:red, 0; green, 0; blue, 0 }  ][line width=0.75]    (10.93,-3.29) .. controls (6.95,-1.4) and (3.31,-0.3) .. (0,0) .. controls (3.31,0.3) and (6.95,1.4) .. (10.93,3.29)   ;
\draw    (265.2,60.4) -- (280.2,60.75) ;
\draw [shift={(282.2,60.8)}, rotate = 181.35] [color={rgb, 255:red, 0; green, 0; blue, 0 }  ][line width=0.75]    (10.93,-3.29) .. controls (6.95,-1.4) and (3.31,-0.3) .. (0,0) .. controls (3.31,0.3) and (6.95,1.4) .. (10.93,3.29)   ;
\draw    (49.2,61.6) -- (92.2,61.79) ;
\draw [shift={(94.2,61.8)}, rotate = 180.25] [color={rgb, 255:red, 0; green, 0; blue, 0 }  ][line width=0.75]    (10.93,-3.29) .. controls (6.95,-1.4) and (3.31,-0.3) .. (0,0) .. controls (3.31,0.3) and (6.95,1.4) .. (10.93,3.29)   ;

\draw (328,53) node [anchor=north west][inner sep=0.75pt]   [align=left] {\begin{minipage}[lt]{54.89pt}\setlength\topsep{0pt}
\begin{center}
{\scriptsize $\dot x\!=\!Ax\!+\!Bu$}
\end{center}

\end{minipage}};
\draw (522,53) node [anchor=north west][inner sep=0.75pt]   [align=left] {\begin{minipage}[lt]{10.09pt}\setlength\topsep{0pt}
\begin{center}
$C$
\end{center}

\end{minipage}};
\draw (386,128) node [anchor=north west][inner sep=0.75pt]   [align=left] {\begin{minipage}[lt]{9.53pt}\setlength\topsep{0pt}
\begin{center}
{\scriptsize $K_{x}$}
\end{center}

\end{minipage}};
\draw (163,46.4) node [anchor=north west][inner sep=0.75pt]    {$\int$};
\draw (219,46) node [anchor=north west][inner sep=0.75pt]   [align=left] {\begin{minipage}[lt]{15.2pt}\setlength\topsep{0pt}
\begin{center}
{\scriptsize $K_{\rm q}$}
\end{center}

\end{minipage}};
\draw (110,75) node [anchor=north west][inner sep=0.75pt]   [align=left] {\begin{minipage}[lt]{8.67pt}\setlength\topsep{0pt}
\begin{center}
$-$
\end{center}

\end{minipage}};
\draw (75,35) node [anchor=north west][inner sep=0.75pt]   [align=left] {\begin{minipage}[lt]{8.68pt}\setlength\topsep{0pt}
\begin{center}
$+$
\end{center}

\end{minipage}};
\draw (260,35) node [anchor=north west][inner sep=0.75pt]   [align=left] {\begin{minipage}[lt]{8.68pt}\setlength\topsep{0pt}
\begin{center}
$+$
\end{center}

\end{minipage}};
\draw (298,75) node [anchor=north west][inner sep=0.75pt]   [align=left] {\begin{minipage}[lt]{8.68pt}\setlength\topsep{0pt}
\begin{center}
$+$
\end{center}

\end{minipage}};
\draw (587,35) node [anchor=north west][inner sep=0.75pt]   [align=left] {\begin{minipage}[lt]{8.67pt}\setlength\topsep{0pt}
\begin{center}
$y$
\end{center}

\end{minipage}};
\draw (470,36) node [anchor=north west][inner sep=0.75pt]   [align=left] {\begin{minipage}[lt]{8.67pt}\setlength\topsep{0pt}
\begin{center}
$x$
\end{center}

\end{minipage}};
\draw (120,35) node [anchor=north west][inner sep=0.75pt]   [align=left] {\begin{minipage}[lt]{8.67pt}\setlength\topsep{0pt}
\begin{center}
$e$
\end{center}

\end{minipage}};
\draw (196,35) node [anchor=north west][inner sep=0.75pt]   [align=left] {\begin{minipage}[lt]{8.67pt}\setlength\topsep{0pt}
\begin{center}
$q$
\end{center}

\end{minipage}};
\draw (308,35) node [anchor=north west][inner sep=0.75pt]   [align=left] {\begin{minipage}[lt]{8.67pt}\setlength\topsep{0pt}
\begin{center}
$u$
\end{center}

\end{minipage}};
\draw (50,35) node [anchor=north west][inner sep=0.75pt]   [align=left] {\begin{minipage}[lt]{11.79pt}\setlength\topsep{0pt}
\begin{center}
$r_{\rm d}$
\end{center}

\end{minipage}};
\end{tikzpicture}
\caption{Full-state feedback control architecture with integral action. }
    \label{fig:FSI_architecture}
\end{figure}

The control law for the full-state feedback controller with integral action (FSI) [p.~550~-~551]\cite{graebe2000control} is
\begin{align}
    v = K_\xi \xi + K_q q,
    \label{eq:v_pp}
\end{align}
where the integrated output error $q(t) \in \BBR^2$ satisfies
\begin{align}
    \dot q 
        =
            r_\rmd - r
        =
            \matl 
                r_{\rmd,1} - r_1 \\
                r_{\rmd,2} - r_2
            \matr.
    \label{eq:v_pp_q}
\end{align}
The gain matrices $K_\xi \in \BBR^{2 \times 8}$ and $K_q \in \BBR^{2 \times 2}$ are chosen using the pole-placement technique to obtain the desired transient response. 
The FSI control architecture is shown in Figure \ref{fig:FSI_architecture}.
The output matrix $C_\rmc$ to compute the output $y=r$ is 
\begin{align}
    C_\rmc 
        \isdef
            \matl 
                1 & 0 & 0& 0 & 0& 0 & 0& 0 \\
                0 & 0 & 0& 0 & 1& 0 & 0& 0 
            \matr.
\end{align}
Note that the FSI control law, given by \eqref{eq:v_pp}, yields asymptotic convergence only when the desired output trajectory consists of step commands. 

\subsubsection{Perfect Tracking Controller (PTC)}

To track a smooth desired trajectory, we consider the control law
%
\begin{align}
    v = K (\xi - \xi_{\rmd}) + B_\rmc ^\rmT \dot \xi_{\rmd}.
    \label{eq:lin_feedback_perfect_tracking}
\end{align}
%
The gain matrix $K \in \BBR^{2 \times 8}$ is chosen using the pole-placement technique to obtain the desired transient response of the state error dynamics 
\begin{align}
    \dot e = (A_\rmc +B_\rmc K) e,
    \label{eq:error_dyn}
\end{align}
where the state error $e \isdef \xi - \xi_\rmd.$ 

Note that the PTC control law  \eqref{eq:lin_feedback_perfect_tracking} requires the desired state trajectory, unlike \eqref{eq:v_pp}, which only requires the desired output trajectory.
However, precisely the additional information encoded in the desired state trajectory allows asymptotically perfect tracking of arbitrary trajectories using the PTC controller  \eqref{eq:lin_feedback_perfect_tracking}.
Furthermore, note that 
\begin{align}
    B_\rmc ^\rmT \dot \xi_\rmd 
        =
            \matl \dot \xi_{\rmd4} \\ \dot \xi_{\rmd8} \matr
        = 
            \matl y_{\rmd1}^{(4)} \\ y_{\rmd2}^{(4)} \matr,
        \label{eq:lin_feed_dot}
\end{align}
which is assumed to be available for smooth trajectories. 
\section{Numerical Simulations}
\label{sec:simulations}

This section presents the numerical simulation results of the bicopter system controlled by the ADIOL controller commanded to track desired trajectories. 
%
In particular, we use the ADIOL controller to track an elliptical trajectory and a second-order Hilbert curve. 

In the following, we assume that the mass and the inertia of the bicopter is $1 \ \rm kg$ and $0.2 \ {\rm kg m^2},$ respectively
%
In the parameter estimator, we set 
$c_1 = 15,$
$c_2 = 5,$
$\alpha_1 = 0.2,$
$\alpha_2 = 1.2,$
$\gamma = 10$, 
$\lambda = 80,$ and 
$\hat \Theta(0) = \matl 2 & 6 \matr^\rmT .$

In the FSI controller, the closed-loop eigenvalues are placed at 
$(-4, -5, -5.5,-6, -8, -4, -5, -5.5,-6, -8)$ 
using Matlab's \texttt{place} routine, which yields
\begin{align}
    K_\xi
        &=
            \matl 
                4876 & 0\\
                1779 & 0\\
                320.5 & 0\\
                28.5 & 0 \\
                0 & 4876\\
                0 & 1779\\
                0 & 320.5\\
                0 & 28.5
            \matr^{\rm T},  
    \
    K_q 
        =
            \matl 
                -528 & 0 \\
                0 & -528
            \matr.
\end{align}

In the PTC controller, the closed-loop eigenvalues are placed at 
$(-1, -2, -3, -5, -1, -2, -3, -5)$
using Matlab's \texttt{place} routine, which yields
\begin{align}
    K 
        &=
            \matl 
                30 & 61 & 41 & 11 & 0 & 0 & 0 & 0\\
                0 & 0 & 0 & 0 & 30 & 61 & 41 & 11
            \matr.
\end{align}





\subsection{Elliptical Trajectory}

\begin{figure}[!ht]
    \centering
    \includegraphics[width=0.8\columnwidth]{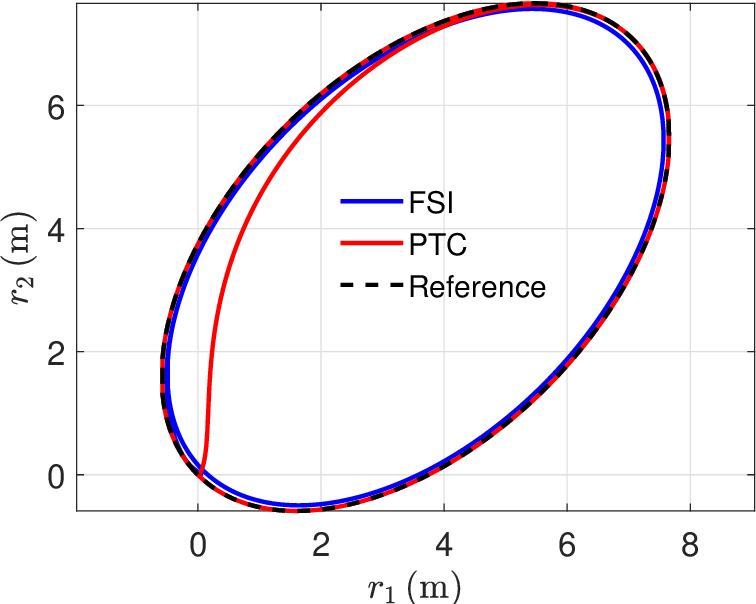}
    \caption{\textbf{Elliptical trajectory}. trajectory-tracking response of the bicopter with ADIOL and two tracking controllers. 
    The output trajectory response with the ADIOL-FSI and ADIOL-PTC controller is shown in blue and red, respectively.
    The reference trajectory is shown in dashed black.}
    \label{fig:IROS23_ADIOL_Bicopter_Trajectory_Circular}
\end{figure}

\begin{figure}[ht]
    \centering
    \includegraphics[width=0.8\columnwidth]{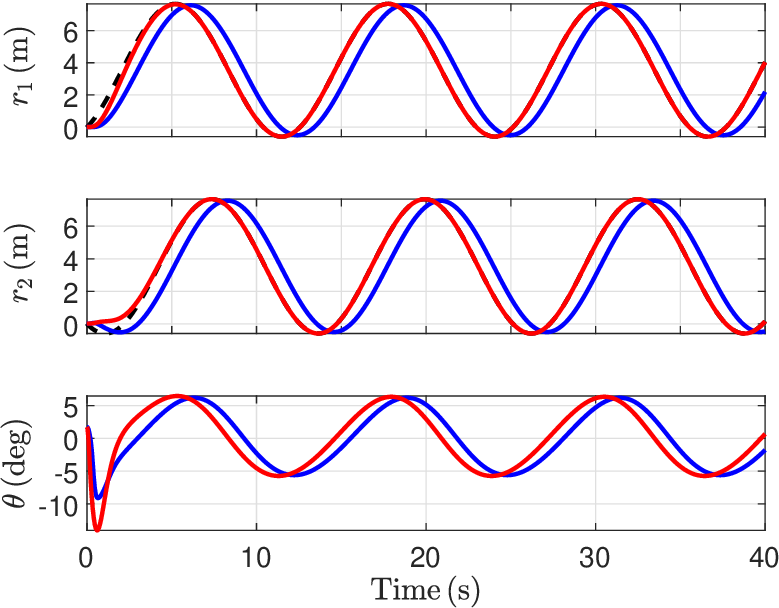}
    \caption{\textbf{Elliptical trajectory}. 
    Position $(r_1, r_2)$ and roll angle $\theta$ response of the bicopter with ADIOL-FSI (in blue) and ADIOL-PTC (in red). 
    }
    \label{fig:IROS23_ADIOL_Bicopter_state}
\end{figure}

\begin{figure}
    \centering
    \includegraphics[width=0.8\columnwidth]{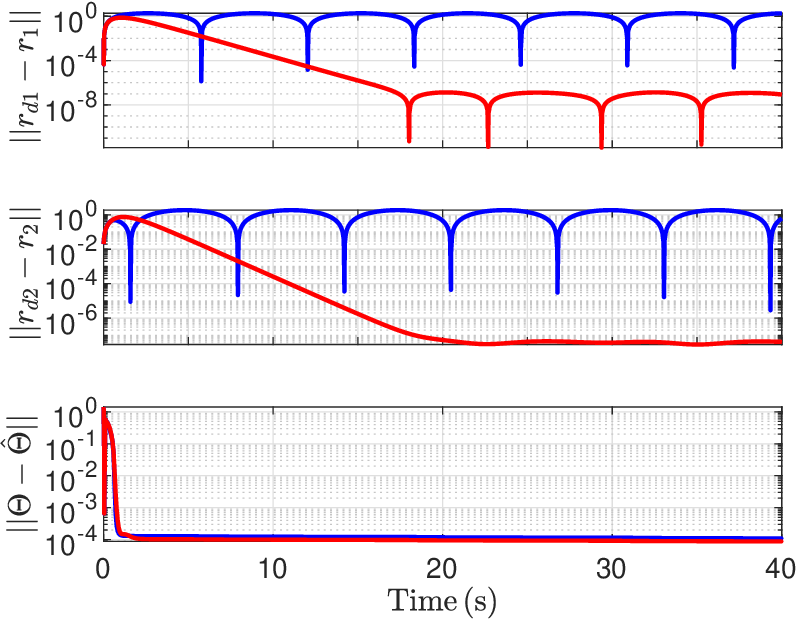}
    \caption{\textbf{Elliptical trajectory}. 
    Position errors and the parameter estimate error with ADIOL-FSI (in blue) and ADIOL-PTC (in red) on a logarithmic scale. 
    }
    \label{fig:IROS23_ADIOL_Bicopter_errors}
\end{figure}

\begin{figure}
    \centering
    \includegraphics[width=0.8\columnwidth]{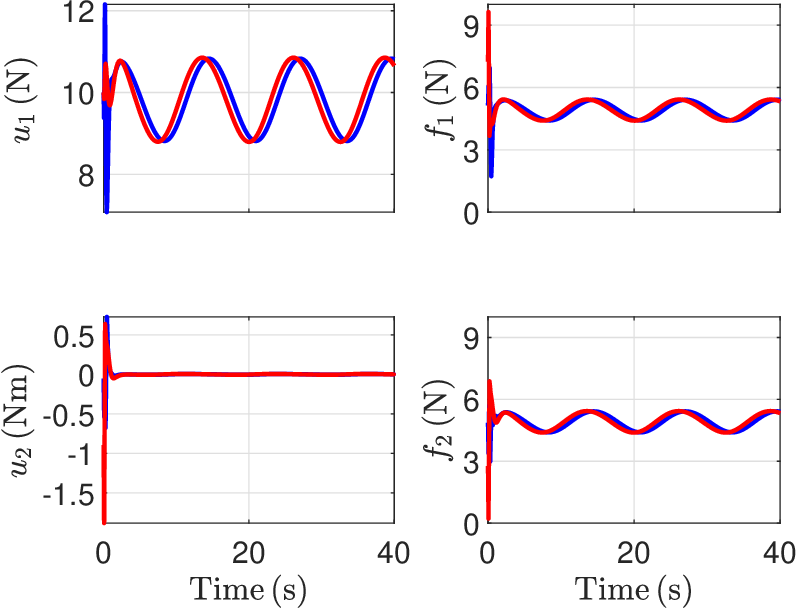}
    \caption{\textbf{Elliptical trajectory}. 
    Control $u$ and the corresponding forces $f_1$ and $f_2$ obtained with ADIOL-FSI (in blue) and ADIOL-PTC (in red).}
    \label{fig:IROS23_ADIOL_Bicopter_control}
\end{figure}

The bicopter is commanded to track a desired elliptical trajectory given by 
\begin{align}
    r_{\rmd1}(t) &= 5 \cos(\phi)-5 \cos(\phi) \cos (\omega t) - 3 \sin(\phi)\sin(\omega t), \nn\\ 
    r_{\rmd2}(t) &= 5\sin(\phi) - 5 \sin(\phi) \cos (\omega t) + 3 \cos(\phi)\sin(\omega t), \nn
\end{align}
where $\phi=45~\rm{deg}$ and $\omega = 0.5 \ \rm rad/s^{-1}. $ 
%
%
Figure \ref{fig:IROS23_ADIOL_Bicopter_Trajectory_Circular} shows the trajectory-tracking response of the bicopter, where the desired trajectory is shown in black dashes and the output trajectory response with the FSI and the PTC tracking controllers is shown in blue and red, respectively.
Figure \ref{fig:IROS23_ADIOL_Bicopter_state} shows the position $r_1, r_2$ response and the roll angle $\theta$ response of the bicopter with both tracking controllers. 
Figure \ref{fig:IROS23_ADIOL_Bicopter_errors} shows the position errors and norm of the parameter estimate error $\Theta-\hat \Theta$ obtained with both tracking controllers on a logarithmic scale. 
Finally, Figure \ref{fig:IROS23_ADIOL_Bicopter_control} shows the control $u$ generated by the ADIOL controller \eqref{eq:IOL_control_law} and the corresponding forces $f_1$ and $f_2.$
Note that the forces $f_1, f_2$ are computed using \eqref{control_signals}.

%

%

\subsection{Second-order Hilbert Curve Trajectory}
\begin{figure}
    \centering
    \includegraphics[width=0.8\columnwidth]{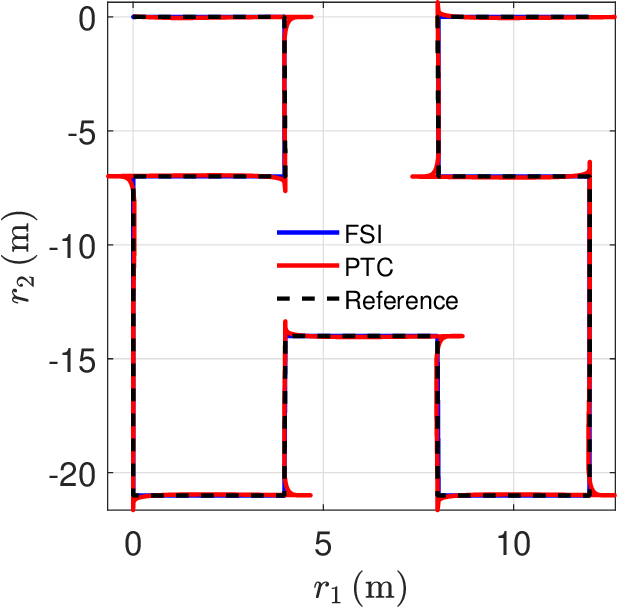}
    \caption{\textbf{Hilbert trajectory}. trajectory-tracking response of the bicopter with ADIOL and two tracking controllers. 
    The output trajectory response with the ADIOL-FSI and ADIOL-PTC controller is shown in blue and red, respectively.
    The reference trajectory is shown in dashed black.}
    \label{fig:IROS23_ADIOL_Bicopter_Trajectory_Hilbert}
\end{figure}

\begin{figure}
    \centering
    \includegraphics[ width=0.8\columnwidth]{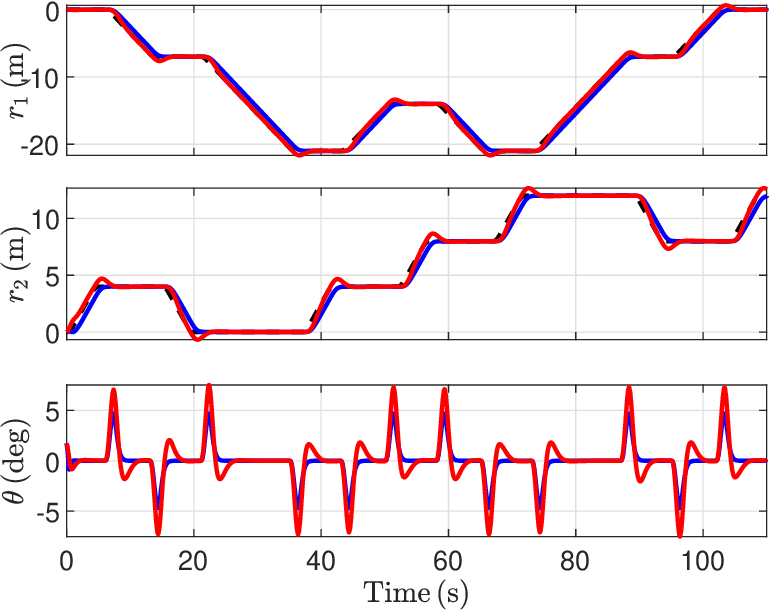}
    \caption{\textbf{Hilbert trajectory}. Position $(r_1, r_2)$ and roll angle $\theta$ response of the bicopter with ADIOL-FSI (in blue) and ADIOL-PTC (in red).}
    \label{fig:IROS23_ADIOL_Bicopter_Hilbert_state}
\end{figure}

\begin{figure}
    \centering
    \includegraphics[width=0.8\columnwidth]{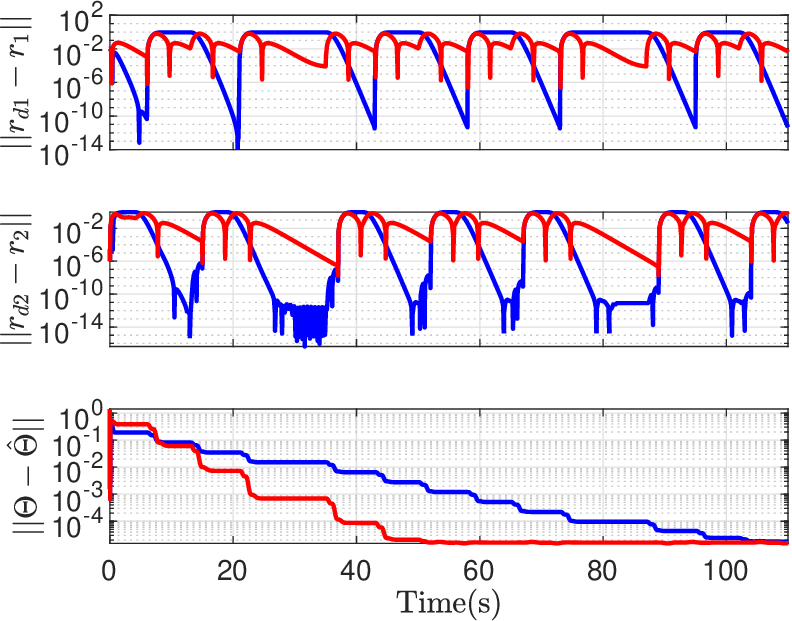}
    \caption{\textbf{Hilbert trajectory}. 
    Position errors and the parameter estimate error with ADIOL-FSI (in blue) and ADIOL-PTC (in red) on a logarithmic scale.}
    \label{fig:IROS23_ADIOL_Bicopter_Hilbert_errors}
\end{figure}

\begin{figure}
    \centering
    \includegraphics[width=0.8\columnwidth]{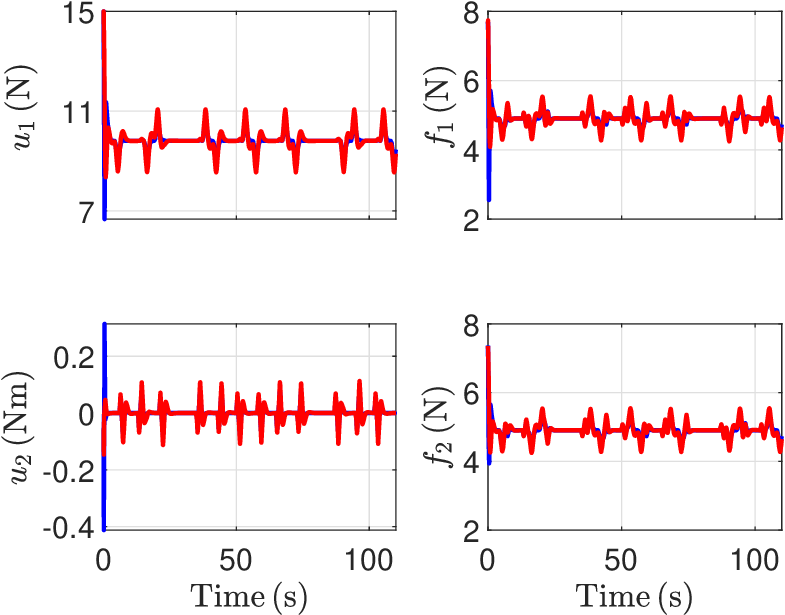}
    \caption{\textbf{Hilbert trajectory}. Control $u$ and the corresponding forces $f_1$ and $f_2$ obtained with ADIOL-FSI (in blue) and ADIOL-PTC (in red).}
    \label{fig:IROS23_ADIOL_Bicopter_Hilbert_control}
    \vspace{-2.0em}
\end{figure}

Next, the bicopter is commanded to track a trajectory constructed using a second-order Hilbert curve.
The trajectory is constructed using the algorithm described in Appendix A of \cite{spencer2022adaptive} with
maximum velocity $v_{\rm max} = 1 \ \rm m/s$ and 
maximum acceleration $a_{\rm max} = 1 \ \rm m/s^2.$ 
%
%
%
%
Figure \ref{fig:IROS23_ADIOL_Bicopter_Trajectory_Hilbert} shows the trajectory-tracking response of the bicopter, where the desired trajectory is shown in black dashes and the output trajectory response with the FSI and the PTC tracking controllers is shown in blue and red, respectively. 
Figure \ref{fig:IROS23_ADIOL_Bicopter_Hilbert_state} shows the position $r_1$ and $r_2$ response and the roll angle $\theta$ response of the bicopter with both tracking contollers. 
Figure \ref{fig:IROS23_ADIOL_Bicopter_Hilbert_errors} shows the position errors, and the norm of the parameter estimate error $\Theta-\hat \Theta$ obtained with both tracking controllers on a logarithmic scale. 
Finally, Figure \ref{fig:IROS23_ADIOL_Bicopter_Hilbert_control} shows the control $u$ generated by the ADIOL controller \eqref{eq:IOL_control_law}, and the corresponding forces $f_1$ and $f_2$. Note that the forces $f_1$ and $f_2$ are computed using \eqref{control_signals}.

The preceding two examples show that the bicopter's trajectory converges to the desired trajectory exponentially without prior perfect knowledge of the bicopter dynamics. 
The parameter estimator monotonically reduces the parameter error. 
However, note that the parameter convergence is not necessary for asymptotic stability of the closed-loop system, as is well known for composite controllers. 



\section{Conclusions and Future Work}
\label{sec:conclusions}
This paper presented an adaptive, dynamic, input-output linearizing controller to linearize the bicopter dynamics and two linear controllers for the bicopter trajectory-tracking problem. 
The controller was constructed by extending the bicopter dynamics to allow linearization, augmenting the controller with a finite-time convergent parameter estimator to yield a composite controller, and designing a linear controller for the linearized dynamics to achieve a desired tracking performance.  
%
%
%
%
The performance of the proposed controller was validated in the numerical simulation of an elliptical and a Hilbert-curve-based trajectory-tracking problem.  

Although the work in this paper does not present a stability analysis of the proposed adaptive controller, numerical simulations confirm that exponentially stable closed-loop dynamics is obtained using the proposed controller.
%
Our future work will focus on relaxing the requirement of full-state feedback by incorporating a finite-time convergent state estimator, characterizing the robustness of the proposed controller to sensor noise, and presenting a rigorous stability analysis of the proposed adaptive control system.

\printbibliography

\end{document}